# ESTIMATION OF A $K$-MONOTONE DENSITY: LIMIT DISTRIBUTION THEORY AND THE SPLINE CONNECTION

By Fadoua Balabdaoui[1] and Jon A. Wellner[2]

*University of Goettingen and University of Washington*

We study the asymptotic behavior of the Maximum Likelihood and Least Squares Estimators of a $k$-monotone density $g_0$ at a fixed point $x_0$ when $k > 2$. We find that the $j$th derivative of the estimators at $x_0$ converges at the rate $n^{-(k-j)/(2k+1)}$ for $j = 0, \ldots, k-1$. The limiting distribution depends on an almost surely uniquely defined stochastic process $H_k$ that stays above (below) the $k$-fold integral of Brownian motion plus a deterministic drift when $k$ is even (odd). Both the MLE and LSE are known to be splines of degree $k-1$ with simple knots. Establishing the order of the random gap $\tau_n^+ - \tau_n^-$, where $\tau_n^\pm$ denote two successive knots, is a key ingredient of the proof of the main results. We show that this "gap problem" can be solved if a conjecture about the upper bound on the error in a particular Hermite interpolation via odd-degree splines holds.

## 1. Introduction.

1.1. *The estimation problem and motivation.* A density function $g$ on $\mathbb{R}^+$ is monotone (or 1-monotone) if it is nonincreasing. It is 2-monotone if it is nonincreasing and convex, and *k-monotone* for $k \geq 3$ if and only if $(-1)^j g^{(j)}$ is nonnegative, nonincreasing and convex for $j = 0, \ldots, k-2$.

We write $\mathcal{D}_k$ for the class of all $k$-monotone densities on $\mathbb{R}^+$ and $\mathcal{M}_k$ for the class of all $k$-monotone functions (without the density restriction). Suppose that $g_0 \in \mathcal{D}_k$ and that $X_1, \ldots, X_n$ are i.i.d. with density $g_0$. We

Received August 2005; revised January 2007.
[1]Supported in part by NSF Grant DMS-02-03320.
[2]Supported in part by NSF Grant DMS-02-03320, NIAID Grant 2R01 AI291968-04 and an NWO grant to the Vrije Universiteit, Amsterdam.

write $\mathbb{G}_n$ for the empirical distribution function of $X_1, \ldots, X_n$. Our main interest is in the Maximum Likelihood Estimators (or MLE's) $\hat{g}_n$ of $g_0 \in \mathcal{D}_k$.

When $k = 1$, it is well known that the maximum likelihood estimator $\hat{g}_n$ of $g_0 \in \mathcal{D}_1$ is the Grenander [14] estimator, that is, the left derivative of the least concave majorant $\widehat{G}_n$ of $\mathbb{G}_n$, and if $g_0'(x_0) < 0$ with $g_0'$ continuous in a neighborhood of $x_0$, then

$$(1.1) \qquad n^{1/3}(\hat{g}_n(x_0) - g_0(x_0)) \xrightarrow{d} (\tfrac{1}{2} g_0(x_0) |g_0'(x_0)|)^{1/3} 2Z,$$

where $2Z$ is the slope at zero of the greatest convex minorant of two-sided Brownian motion $+ t^2$, $t \in \mathbb{R}$; see Prakasa Rao [35], Groeneboom [15] and Kim and Pollard [24].

When $k = 2$, Groeneboom, Jongbloed and Wellner [18] considered both the MLE and LSE and established that if the true convex and nonincreasing density $g_0$ satisfies $g_0''(x_0) > 0$ (and $g_0''$ is continuous in a neighborhood of $x_0$), then

$$(1.2) \quad \begin{pmatrix} n^{2/5}(\bar{g}_n(x_0) - g_0(x_0)) \\ n^{1/5}(\bar{g}_n'(x_0) - g'(x_0)) \end{pmatrix} \xrightarrow{d} \begin{pmatrix} (\tfrac{1}{24} g_0^2(x_0) g_0''(x_0))^{1/5} H^{(2)}(0) \\ (\tfrac{1}{24^3} g_0(x_0) g_0''(x_0)^3)^{1/5} H^{(3)}(0) \end{pmatrix},$$

where $\bar{g}_n$ is either the MLE or LSE and $H$ is a random cubic spline function such that $H^{(2)}$ is convex and $H$ stays above integrated two-sided Brownian motion $+ t^4$, $t \in \mathbb{R}$, and touches exactly at those points where $H^{(2)}$ changes its slope; see Groeneboom, Jongbloed and Wellner [17].

Our main interest in this paper is in establishing a generalization of the pointwise limit theory given in (1.1) and (1.2) for general $k \in \mathbb{N}$, $k \geq 1$.

Beyond the obvious motivation of extending the known results for $k = 1$ and $k = 2$ as listed above, there are several further reasons for considering such extensions:

(a) Pointwise limit distribution theory for natural nonparametric estimators of the piecewise smooth regression models of smoothness $k$ considered by Mammen [29] is only available for $k \in \{1, 2\}$. Similar models (with just one element in the partition) have been proposed for software reliability problems by Miller and Sofer [33]. Similarly, pointwise limit distribution theory is still lacking for the locally adaptive regression spline estimators considered by Mammen and van de Geer [30].

(b) The classes of densities $\mathcal{D}_k$ have mixture representations as scale mixtures of Beta$(1, k)$ densities: as is known from Williamson [43] (see also Lévy [26], Gneiting [13] and Balabdaoui and Wellner [2]), $g \in \mathcal{D}_k$ if and only if there is a distribution function $F$ on $(0, \infty)$ such that

$$(1.3) \qquad g(x) = \int_0^\infty \frac{k}{y^k}(y - x)_+^{k-1} dF(y) = \int_0^\infty w \left(1 - \frac{wx}{k}\right)_+^{k-1} d\tilde{F}(w),$$

where $z_+ \equiv z1\{z \geq 0\}$ and $\tilde{F} = F(k/\cdot)$. The second form of the mixture representation in the last display makes it clear that the limiting class of densities as $k \to \infty$, namely $\mathcal{D}_\infty$, is the class of scale mixtures of exponential distributions. In view of Feller [11], pages 232–233, this is just the class of *completely monotone* densities; see also Widder [42] and Gneiting [12]. To the best of our knowledge, there is no pointwise limit distribution theory available for the MLE in any class of mixed densities based on a smooth mixing kernel, including this particular case in which the kernel (or mixture density) is the exponential scale family as studied by Jewell [22]. On the other hand, maximum likelihood estimators in various classes of mixture models with smooth kernels have been proposed in a wide range of applications including pharmacokinetics (Mallet [27], Mallet, Mentré, Steimer and Lokiec [28] and Davidian and Gallant [6]), demography (Vaupel, Manton and Stallard [41]) and shock models and variations in hazard rates (Harris and Singpurwalla [20], McNolty, Doyle and Hansen [31] and Hill, Saunders and Laud [21]).

(c) The whole family of mixture models $\mathcal{D}_k$ corresponding to $k \in (0, \infty)$ in (1.3) might eventually be of some interest, especially since the family of distributions corresponding to the classical Wicksell problem is contained in the class $\mathcal{D}_{1/2}$; see, for example, Groeneboom and Jongbloed [16].

(d) The subclass of $k$-monotone densities with mixing distribution $F$ satisfying $g^{(k-1)}(0) = k! \int_0^\infty y^{-k} \, dF(y) < \infty$ can be regarded as the class of distributions arising in a generalization of Hampel's bird-watching problem (Hampel [19]), in which birds are captured $k$ times, but only one "intercatch" time is recorded. Based on those observed intercatch times, the goal is to estimate the true distribution $F$ of the resting times $Y$ of the migrating birds, which we assume to have a density $f$ with $k$th moment $\mu_k(f) < \infty$. Furthermore, we assume that the time points of capture form the arrival time points of a Poisson process with rate $\lambda$, that given $Y = y$, the number of captures by time $y$ is Poisson$(\lambda y)$ with $\lambda$ small enough so that $\exp(-\lambda y) \approx 1$ and that the probability of catching a bird more than $k$ times is negligible (see also Hampel [19] and Anevski [1]). If $S_{k,1}$ denotes the elapsed time between the first and second captures (the only observed intercatch time), then it follows by a derivation analogous to Hampel's that the density of the time $S_{k,1}$ is given by

$$g(x) = \frac{1}{\mu_k(f)} \int_0^\infty k(y - x)_+^{k-1} f(y) \, dy,$$

which is clearly $k$-monotone. We obtain $F$, the probability distribution of $Y$, by inverting the previous mixture representation, that is,

$$F(t) = 1 - \frac{g^{(k-1)}(t)}{g^{(k-1)}(0+)}$$

at any point of continuity $t > 0$ of $F$.

In connection with (a), it is interesting to note that the definition of the family $\mathcal{D}_k$ is equivalent to $g \in \mathcal{D}_k$ if and only if $(-1)^{k-1} g^{(k-1)}$ (where $g^{(k-1)}$ is either the left or right derivative of $g^{(k-2)}$) is nonincreasing. This follows from Lemma 4.3 of Gneiting [13] since Gneiting's condition $\lim_{x \to \infty} g(x) = 0$ is automatic for densities. Thus the equivalent definition of $\mathcal{D}_k$ has a natural connection with the work of Mammen [29] in the nonparametric regression setting. In parallel to the treatment of convex regression estimation given by Groeneboom, Jongbloed and Wellner [18], it seems clear that pointwise distribution theory for nonparametric least squares estimators for the regression problems in (a) could be developed if adequate theory were available for the Maximum Likelihood and Least Squares estimators of densities in the class $\mathcal{D}_k$, so we focus exclusively on the density case in this paper. In Section 5, we comment further on the difficulties in obtaining corresponding limit theory for the smooth kernel cases discussed in (b).

1.2. *Description of the key difficulty*: *the gap problem.* The key result that Groeneboom, Jongbloed and Wellner [18] used to establish (1.2) is that $\tau_n^+ - \tau_n^- = O_p(n^{-1/5})$ as $n \to \infty$, where $\tau_n^-$ and $\tau_n^+$ are two successive jump points of the first derivative of $\bar{g}_n$ in the neighborhood of $x_0$. Such a result was already proved by Mammen [29] (see Lemma 8) in the context of nonparametric regression, where the true regression curve, $m$, is piecewise concave/convex or convex/concave such that $m$ is twice continuously differentiable in the neighborhood of $x_0$, and $m''(x_0) \neq 0$. Furthermore, Mammen [29] conjectured the right form of the asymptotic distribution of his Least Squares estimator, which was later established by Groeneboom, Jongbloed and Wellner [18].

To obtain the stochastic order $n^{-1/5}$ for the gap, Groeneboom, Jongbloed and Wellner [18] used the characterizations of the estimators, together with the "midpoint property" which we review in Section 4. For $k = 1$, the same property can be used to establish that $n^{-1/3}$ is the order of the gap. As a function of $k$, it is natural to conjecture that $n^{-1/(2k+1)}$ is the general form of the order of the gap. In the problem of nonparametric regression via splines, Mammen and van de Geer [30] have conjectured that $n^{-1/(2k+1)}$ is the order of the distance between the knot points of their regression spline $\hat{m}$ under the assumption that the true regression curve $m_0$ satisfies our same working assumptions, but the question was left open (see Mammen and van de Geer [30], page 400). In this paper, we refer to the problem of establishing the order of $\tau_n^+ - \tau_n^-$ as the *gap problem*.

In Section 4, we show that when $k > 2$, the gap problem is closely related to a "nonclassical" Hermite interpolation problem via odd-degree splines. To put the interpolation problem encountered in the next section in context, it is useful to review briefly the related *complete interpolation problem* for odd-degree splines which is more "classical" and for which error bounds uniform in the knots are now available. Given a function $f \in C^{(k-1)}[0,1]$ and an increasing sequence $0 = y_0 < y_1 < \cdots < y_m < y_{m+1} = 1$, where $m \geq 1$ is an integer, it is well known that there exists a unique spline, called the *complete spline* and denoted here by $Cf$, of degree $2k-1$ with interior knots $y_1, \ldots, y_m$ that satisfies the $2k + m$ conditions

$$(Cf)(y_i) = f(y_i), \qquad i = 1, \ldots, m,$$
$$(Cf)^{(l)}(y_0) = f^{(l)}(y_0), \qquad (Cf)^{(l)}(y_{m+1}) = f^{(l)}(y_{m+1}), \qquad l = 0, \ldots, k-1;$$

see Schoenberg [36], de Boor [8] or Nürnberger [34], page 116, for further discussion. If $j \in \{0, \ldots, k\}$ and $f \in C^{(k+j)}[0,1]$, then there exists $c_{k,j} > 0$ such that

$$(1.4) \qquad \sup_{0 < y_1 < \cdots < y_m < 1} \|f - Cf\|_\infty \leq c_{k,j} \|f^{(k+j)}\|_\infty.$$

For $j = k$, this "uniform in knots" bound in the complete interpolation problem was first conjectured by de Boor [7] for $k > 4$ as a generalization that goes beyond $k = 2, 3$ and $4$, for which the result was already established (see also de Boor [8]). By a scaling argument, the bound (1.4) implies that if $f \in C^{(2k)}[a,b], a < b \in \mathbb{R}$, then the interpolation error in the complete interpolation problem is uniformly bounded in the knots and the bound is of the order of $(b-a)^{2k}$. One key property of the complete spline interpolant $Cf$ is that $(Cf)^{(k)}$ is the Least Squares approximation of $f^{(k)}$ when $f^{(k)} \in L_2([0,1])$, that is, if $\mathcal{S}_k(y_1, \ldots, y_m)$ denotes the space of splines of order $k$ (degree $k-1$) and interior knots $y_1, \ldots, y_m$, then

$$(1.5) \quad \int_0^1 ((Cf)^{(k)} - f^{(k)}(x))^2 \, dx = \min_{S \in \mathcal{S}_k(y_1, \ldots, y_m)} \int_0^1 (S(x) - f^{(k)}(x))^2 \, dx$$

(see, e.g., Schoenberg [36], de Boor [8], Nürnberger [34]). Consequently, if $L_\infty$ denotes the space of bounded functions on $[0,1]$, then the properly defined map

$$C^{(k)}[0,1] \to \mathcal{S}_k(\underline{y}),$$
$$f^{(k)} \to (Cf)^{(k)},$$

where $\underline{y} = (y_1, \ldots, y_m)$, is the restriction of the orthoprojector $P_{S_k(\underline{y})}$ from $L_\infty$ to $\mathcal{S}_k(\underline{y})$ with respect to the inner product $\langle g, h \rangle = \int_0^1 g(x) h(x) \, dx$ which assigns to a function $g \in L_\infty$ the $k$-*th* derivative of the complete spline interpolant of *any* primitive of $g$ of order $k$ (note that the difference between two primitives of $g$ of order $k$ is a polynomial of degree $k-1$).

de Boor [8] pointed out that in order to prove the conjecture, it is enough to prove that

$$\sup_{\underline{y}} \|P_{S_k(\underline{y})}\|_\infty = \sup_{\underline{y}} \sup_{g \in L_\infty} \frac{\|P_{S_k(\underline{y})}(g)\|_\infty}{\|g\|_\infty}$$

is bounded. This was successfully achieved by Shadrin [38].

The Hermite interpolation problem which arises naturally in Section 4 appears to be another variant of interpolation problems via odd-degree splines which has not yet been studied in the approximation theory or spline literature. More specifically, if $f$ is some real-valued function in $C^{(j)}[0,1]$ for some $j \geq 1$ and $0 = y_0 < y_1 < \cdots < y_{2k-4} < y_{2k-3} = 1$ is a given increasing sequence, then there exist a unique spline $\mathcal{H}_k f$ of degree $2k-1$ and interior knots $y_1, \ldots, y_{2k-4}$ satisfying the $4k-4$ conditions

(1.6) $(\mathcal{H}_k f)(y_i) = f(y_i)$ and $(\mathcal{H}_k f)'(y_i) = f'(y_i)$, $i = 0, \ldots, 2k-3$.

It turns out that deriving the stochastic order of the distance between two successive knots of the MLE and LSE in the neighborhood of the point of estimation is very closely linked to bounding the error in this new Hermite interpolation independently of the locations of the knots of the spline interpolant. More precisely, if $g_t(x) = (x-t)_+^{k-1}/(k-1)!$ is the power truncated function of degree $k-1$ with unique knot $t$, then we conjecture that there is a constant $d_k > 0$ such that

(1.7) $$\sup_{t \in (0,1)} \sup_{0 < y_1 < \cdots < y_{2k-4} < 1} \|g_t - \mathcal{H}_k g_t\|_\infty \leq d_k.$$

As shown in Balabdaoui and Wellner [3], the preceding formulation implies that boundedness of the error uniformly in the knots of the spline interpolant holds true for any $f \in C^{(k+j)}$, that is,

$$\sup_{0 < y_1 < \cdots < y_{2k-4} < 1} \|f - \mathcal{H}_k f\|_\infty \leq d_{k,j} \|f^{(k+j)}\|_\infty.$$

If $j = k$ and $\|f^{(2k)}\|_\infty \leq 1$, it follows from Proposition 1 of Balabdaoui and Wellner [3] that the interpolation error must be bounded above by the error for interpolating the perfect spline,

$$S^*(t) = \frac{1}{(2k)!}\left(t^{2k} + 2\sum_{i=1}^{2k-4} (-1)^i (t - \tau_j)_+^{2k}\right).$$

For a definition of perfect splines, see, for example, Bojanov, Hakopian and Sahakian [5], Chapter 6. Based on a large number of simulations, we found that

$$\sup_{0 < y_1 < \cdots < y_{2k-4} < 1} \|S^* - \mathcal{H}_k S^*\|_\infty \leq \frac{2}{(2k)!}$$

for fairly large values of $k$ (see the last column in Table 2 in Balabdaoui and Wellner [3]). The latter strongly suggests that for $f \in C^{(2k)}[0,1]$, we have

$$(1.8) \qquad \sup_{0<y_1<\cdots<y_{2k-4}<1} \|f - \mathcal{H}_k f\|_\infty \leq \frac{2}{(2k)!} \|f^{(2k)}\|_\infty.$$

Based on conjecture (1.7), we will prove that the distance between two consecutive knots in a neighborhood of $x_0$ is $O_p(n^{-1/(2k+1)})$.

After a brief introduction to the MLE and LSE and their respective characterizations, we give in Section 3 a statement of our main result which gives the joint asymptotic distribution of the successive derivatives of the MLE and LSE. The obtained convergence rate $n^{-(k-j)/(2k+1)}$ for the $j$th derivative of any of the estimators was found by Balabdaoui and Wellner [2] to be the asymptotic minimax lower bound for estimating $g_0^{(j)}(x_0)$, $j = 0, \ldots, k-1$, under the same working assumptions. The limiting distribution depends on the higher derivatives of $H_k$, an almost surely uniquely defined process that stays above (below) the $(k-1)$-fold integral of Brownian motion plus the drift $(k!/(2k)!)t^{2k}$ when $k$ is even (odd) and whose derivative of order $2k-2$ is convex [$H_k$ is also said to be $(2k-2)$-*convex*]. The process $H_k$ is studied separately in Balabdaoui and Wellner [2]. Proving the existence of $H_k$ also relies on our conjecture in (1.7) since the key problem, also referred to as the *gap problem*, depends on a very similar Hermite interpolation problem, except that the knots of the estimators are replaced by the points of touch between the $(k-1)$-fold integral of Brownian motion plus the drift $(k!/(2k)!)t^{2k}$ and $H_k$. For more discussion of the background and related problems, see Balabdaoui and Wellner [2]. For a discussion of algorithms and computational issues, see Balabdaoui and Wellner [2].

**2. The estimators and their characterization.** Let $X_1, \ldots, X_n$ be $n$ independent observations from a common $k$-monotone density $g_0$. We consider nonparametric estimation of $g_0$ via the Least Squares and Maximum Likelihood methods, and that of its mixture distribution $F_0$, that is, the distribution function on $(0, \infty)$ such that

$$g_0(x) = \int_0^\infty \frac{k(t-x)_+^{k-1}}{t^k} \, dF_0(t), \qquad x > 0.$$

In other words, $g_0$ is a scale mixture of Beta$(1,k)$ densities. The mixing distribution is, furthermore, given at any point of continuity $t$ by the inversion formula

$$(2.1) \qquad F_0(t) = \sum_{j=0}^k (-1)^j \frac{t^j}{j!} G_0^{(j)}(t),$$

where $G_0(t) = \int_0^t g_0(x)\,dx$. An estimator for $F_0$ can be obtained by simply plugging in estimators of $G_0^{(j)} = g_0^{(j-1)}$, $j = 0, \ldots, k$, in the inversion formula (2.1). We call estimation of the (mixed) $k$-monotone density $g_0$ the *direct problem* and estimation of the mixing distribution function $F_0$ the *inverse problem*. For more technical details on the mixture representation and the inversion formula, see Lemma 2.1 of Balabdaoui and Wellner [2].

We now give the definitions of the Least Squares and Maximum Likelihood estimators; these were already considered in the case $k = 2$ by Groeneboom, Jongbloed and Wellner [18]. The LSE, $\tilde{g}_n$, is the minimizer of the criterion function

$$\Phi_n(g) = \frac{1}{2}\int_0^\infty g^2(t)\,dt - \int_0^\infty g(t)\,d\mathbb{G}_n(t)$$

over the class $\mathcal{M}_k$, whereas the MLE, $\hat{g}_n$, maximizes the "adjusted" log-likelihood function, that is,

$$l_n(g) = \int_0^\infty \log g(t)\,d\mathbb{G}_n(t) - \int_0^\infty g(t)\,dt,$$

over the same class. In Balabdaoui and Wellner [2], we find that both estimators exist and are splines of degree $k-1$, that is, their $(k-1)$st derivative is stepwise. Furthermore, as shown in Balabdaoui and Wellner [2], the LSE's and MLE's are characterized as follows: let $\tilde{H}_n$ and $\mathbb{Y}_n$ be the processes defined for all $x \geq 0$ by

$$\begin{aligned}(2.2)\quad \mathbb{Y}_n(x) &= \int_0^x \int_0^{t_{k-1}} \cdots \int_0^{t_2} \mathbb{G}_n(t_1)\,dt_1\,dt_2\cdots dt_{k-1} \\ &= \int_0^x \frac{(x-t)^{k-1}}{(k-1)!}\,d\mathbb{G}_n(t)\end{aligned}$$

and

$$\begin{aligned}(2.3)\quad \tilde{H}_n(x) &= \int_0^x \int_0^{t_k} \cdots \int_0^{t_2} \tilde{g}_n(t_1)\,dt_1\,dt_2\cdots dt_k \\ &= \int_0^x \frac{(x-t)^{k-1}}{(k-1)!} \tilde{g}_n(t)\,dt.\end{aligned}$$

Then the $k$-monotone function $\tilde{g}_n$ is the LSE if and only if

$$(2.4)\quad \tilde{H}_n(x) \begin{cases} \geq \mathbb{Y}_n(x), & \text{for all } x \geq 0, \\ = \mathbb{Y}_n(x), & \text{if } (-1)^{k-1}\tilde{g}_n^{(k-1)}(x-) < (-1)^{k-1}\tilde{g}_n^{(k-1)}(x+). \end{cases}$$

For the MLE, we define the process

$$(2.5)\quad \widehat{H}_n(x, g) = \int_0^x \frac{k(x-t)^{k-1}}{x^k \hat{g}_n(t)}\,d\mathbb{G}_n(t)$$

for all $x \geq 0$ and $g \in \mathcal{D}_k$. A necessary and sufficient condition for the $k$-monotone function $\hat{g}_n$ to be the MLE is then given by

$$(2.6) \quad \widehat{H}_n(x, \hat{g}_n) \begin{cases} \leq 1, & \text{for all } x \geq 0, \\ = 1, & \text{if } (-1)^{k-1}\hat{g}_n^{(k-1)}(x-) < (-1)^{k-1}\hat{g}_n^{(k-1)}(x+). \end{cases}$$

These characterizations are crucial for understanding the local asymptotic behavior of the LSE and MLE. They were exploited in Balabdaoui and Wellner [2] to show uniform strong consistency of the estimators on intervals of the form $[c, \infty), c > 0$. Here, they prove to be once again very useful for establishing the limit theory in both the direct and inverse problems.

### 3. The asymptotic distribution.

3.1. *The main convergence theorem.* To prepare for a statement of the main result, we first recall the following theorem from Balabdaoui and Wellner [2] giving existence of the processes $H_k$.

THEOREM 3.1. *For all $k \geq 1$, let $Y_k$ denote the stochastic process defined by*

$$Y_k(t) = \begin{cases} \int_0^t \dfrac{(t-s)^{k-1}}{(k-1)!} \, dW(s) + \dfrac{(-1)^k k!}{(2k)!} t^{2k}, & t \geq 0, \\ \int_t^0 \dfrac{(t-s)^{k-1}}{(k-1)!} \, dW(s) + \dfrac{(-1)^k k!}{(2k)!} t^{2k}, & t < 0. \end{cases}$$

*If conjecture (1.7) holds (see also the discussion in Balabdaoui and Wellner [2]), then there exists an almost surely uniquely defined stochastic process $H_k$ characterized by the following four conditions:*

(i) *the process $H_k$ stays everywhere above the process $Y_k$:*

$$H_k(t) \geq Y_k(t), \qquad t \in \mathbb{R};$$

(ii) $(-1)^k H_k$ *is $2k$-convex, that is, $(-1)^k H_k^{(2k-2)}$ exists and is convex;*
(iii) *the process $H_k$ satisfies*

$$\int_{-\infty}^{\infty} (H_k(t) - Y_k(t)) \, dH_k^{(2k-1)}(t) = 0;$$

(iv) *if $k$ is even, $\lim_{|t|\to\infty}(H_k^{(2j)}(t) - Y_k^{(2j)}(t)) = 0$ for $j = 0, \ldots, (k-2)/2$; if $k$ is odd, $\lim_{t\to\infty}(H_k(t) - Y_k(t)) = 0$ and $\lim_{|t|\to\infty}(H_k^{(2j+1)}(t) - Y_k^{(2j+1)}(t)) = 0$ for $j = 0, \ldots, (k-3)/2$.*

We are now able to state the main result of this paper, which generalizes Theorem 6.2 of Groeneboom, Jongbloed and Wellner [18] for estimating convex (2-monotone) densities.

THEOREM 3.2. *Let $x_0 > 0$ and $g_0$ be a $k$-monotone density such that $g_0$ is $k$-times differentiable at $x_0$ with $(-1)^k g_0^{(k)}(x_0) > 0$ and assume that $g_0^{(k)}$ is continuous in a neighborhood of $x_0$. Let $\bar{g}_n$ denote either the LSE $\tilde{g}_n$ or the MLE $\hat{g}_n$ and let $\bar{F}_n$ be the corresponding mixing measure defined in terms of $\bar{G}_n = \int_0^\cdot \bar{g}_n(s)\,ds$ via (2.1). If conjecture (1.7) holds, then*

$$\begin{pmatrix} n^{k/(2k+1)}(\bar{g}_n(x_0) - g_0(x_0)) \\ n^{(k-1)/(2k+1)}(\bar{g}_n^{(1)}(x_0) - g_0^{(1)}(x_0)) \\ \vdots \\ n^{1/(2k+1)}(\bar{g}_n^{(k-1)}(x_0) - g_0^{(k-1)}(x_0)) \end{pmatrix} \xrightarrow{d} \begin{pmatrix} c_0(x_0) H_k^{(k)}(0) \\ c_1(x_0) H_k^{(k+1)}(0) \\ \vdots \\ c_{k-1}(x_0) H_k^{(2k-1)}(0) \end{pmatrix}$$

*and*

$$n^{1/(2k+1)}(\bar{F}_n(x_0) - F(x_0)) \xrightarrow{d} \frac{(-1)^k x_0^k}{k!} c_{k-1}(x_0) H_k^{(2k-1)}(0),$$

*where*

$$c_j(x_0) = \left\{ (g_0(x_0))^{k-j} \left( \frac{(-1)^k g_0^{(k)}(x_0)}{k!} \right)^{2j+1} \right\}^{1/(2k+1)},$$

*for $j = 0, \ldots, k-1$.*

3.2. *The key results and outline of the proofs.* Our proof of Theorem 3.2 proceeds by solving the key gap problem assuming that our conjecture (1.7) holds. This is carried out in Section 4 in which the main result is the following.

LEMMA 3.1. *Let $k \geq 3$ and $\bar{g}_n$ denote either the LSE $\tilde{g}_n$ or the MLE $\hat{g}_n$. If $g_0 \in \mathcal{D}_k$ satisfies $g_0^{(k)}(x_0) \neq 0$ and conjecture (1.7) holds, then $\tau_{2k-3} - \tau_0 = O_p(n^{-1/(2k+1)})$, where $\tau_0 < \cdots < \tau_{2k-3}$ are $2k-2$ successive jump points of $\bar{g}_n^{(k-1)}$ in a neighborhood of $x_0$.*

Using Lemma 3.1, we can establish the rate(s) of convergence of the estimators $\tilde{g}_n$ and $\hat{g}_n$ and their derivatives viewed as local processes in $n^{-1/(2k+1)}$ neighborhoods of the fixed point $x_0$. This is accomplished in Proposition 3.1. Once the rates have been established, we define for the LSE localized versions $\mathbb{Y}_n^{\text{loc}}$, $\tilde{H}_n^{\text{loc}}$ of the processes $\mathbb{Y}_n$, $\tilde{H}_n$ given in (2.2) and (2.3), respectively, and $\widehat{\mathbb{Y}}_n^{\text{loc}}$, $\widehat{H}_n^{\text{loc}}$ related to the process $\widehat{H}_n$ given in (2.5) in the case of the MLE. The proof then proceeds by showing that:

- the localized processes $\mathbb{Y}_n^{\text{loc}}$ and $\widehat{\mathbb{Y}}_n^{\text{loc}}$ converge weakly to $\mathbb{Y}_{a,\sigma}$, where

$$Y_{a,\sigma}(t) = \begin{cases} \sigma \int_0^t \int_0^{s_{k-1}} \cdots \int_0^{s_2} W(s_1) \, ds_1 \cdots ds_{k-1} + a(-1)^k \dfrac{k!}{(2k)!} t^{2k}, \\ \quad t \geq 0, \\ \sigma \int_t^0 \int_{s_{k-1}}^0 \cdots \int_{s_2}^0 W(s_1) \, ds_1 \cdots ds_{k-1} + a(-1)^k \dfrac{k!}{(2k)!} t^{2k}, \\ \quad t \leq 0, \end{cases}$$

with $\sigma = \sqrt{g(x_0)}$, $a = (-1)^k g_0^{(k)}(x_0)/k!$, and where $W$ is a two-sided Brownian motion process starting from 0; this can be shown by classical methods from Shorack and Wellner [39] or alternatively via the strong approximation of Komlós, Major and Tusnády [25];

- the localized processes $\tilde{H}_n^{\text{loc}}$ and $\widehat{H}_n^{\text{loc}}$ satisfy Fenchel (inequality and equality) relations relative to the localized processes $\mathbb{Y}_n^{\text{loc}}$ and $\widehat{\mathbb{Y}}_n^{\text{loc}}$, respectively.

We then show via tightness that the localized processes $\tilde{H}_n^{\text{loc}}$ and $\widehat{H}_n^{\text{loc}}$ (and all their derivatives up to order $2k-1$) converge to a limit process satisfying the conditions (i)–(iv) of Theorem 3.1 and hence the limit process in both cases is just $H_k$ (up to scaling by constants). When specialized to $t=0$, this gives the conclusion of Theorem 3.2.

The following is the key proposition concerning rates of convergence.

PROPOSITION 3.1. *Fix $x_0 > 0$ and let $g_0$ be a $k$-monotone density such that $(-1)^k g_0^{(k)}(x_0) > 0$. Let $\bar{g}_n$ denote either the MLE $\hat{g}_n$ or the LSE $\tilde{g}_n$. If conjecture (1.7) holds, then for each $M > 0$, we have*

$$\sup_{|t| \leq M} \left| \bar{g}_n^{(j)}(x_0 + n^{-1/(2k+1)} t) - \sum_{i=j}^{k-1} \frac{n^{-(i-j)/(2k+1)} g_0^{(i)}(x_0)}{(i-j)!} t^{i-j} \right|$$

(3.1)

$$= O_p(n^{-(k-j)/(2k+1)}) \quad \text{for } j = 0, \ldots, k-1.$$

For the LSE, we define the local $\mathbb{Y}_n$- and $\tilde{H}_n$-processes by

$$\mathbb{Y}_n^{\text{loc}}(t) = n^{2k/(2k+1)} \int_{x_0}^{x_0 + tn^{-1/(2k+1)}} \int_{x_0}^{v_{k-1}} \cdots \int_{x_0}^{v_2} \bigg\{ \mathbb{G}_n(v_1) - \mathbb{G}_n(x_0)$$

$$- \int_{x_0}^{v_1} \sum_{j=0}^{k-1} \frac{(u - x_0)^j}{j!} g_0^{(j)}(x_0) \, du \bigg\} \prod_{i=1}^{k-1} dv_i$$

and

$$\tilde{H}_n^{\text{loc}}(t) = n^{2k/(2k+1)} \int_{x_0}^{x_0 + tn^{-1/(2k+1)}} \int_{x_0}^{v_k} \cdots \int_{x_0}^{v_2} \bigg\{ \tilde{g}_n(v_1)$$

$$-\sum_{j=0}^{k-1} \frac{(v_1-x_0)^j}{j!} g_0^{(j)}(x_0) \bigg\} dv_1 \cdots dv_k$$

$$+ \tilde{A}_{k-1,n} t^{k-1} + \tilde{A}_{k-2,n} t^{k-2} + \cdots + \tilde{A}_{1,n} t + \tilde{A}_{0,n},$$

respectively, where

$$\tilde{A}_{j,n} = \frac{n^{(2k-j)/(2k+1)}}{j!}(\tilde{H}_n^{(j)}(x_0) - \mathbb{Y}_n^{(j)}(x_0)), \qquad j = 0, \ldots, k-1.$$

Let $r_k \equiv 1/(2k+1)$. In the case of the MLE, the local processes $\widehat{\mathbb{Y}}_n^{\text{loc}}$ and $\widehat{H}_n^{\text{loc}}$ are defined as

$$\frac{\widehat{\mathbb{Y}}_n^{\text{loc}}(t)}{g_0(x_0)} = n^{2kr_k} \int_{x_0}^{x_0+tn^{-r_k}} \int_{x_0}^{v_{k-1}} \cdots \int_{x_0}^{v_1} \frac{g_0(v) - \sum_{j=0}^{k-1}(v-x_0)^j/j! g_0^{(j)}(x_0)}{\hat{g}_n(v)}$$
$$dv\, dv_1 \cdots dv_{k-1}$$
$$+ n^{2kr_k} \int_{x_0}^{x_0+tn^{-r_k}} \int_{x_0}^{v_{k-1}} \cdots \int_{x_0}^{v_1} \frac{1}{\hat{g}_n(v)} d(\mathbb{G}_n - G_0)(v) dv_1 \cdots dv_{k-1}$$

and

$$\frac{\widehat{H}_n^{\text{loc}}(t)}{g_0(x_0)} = n^{2kr_k} \int_{x_0}^{x_0+tn^{-r_k}} \int_{x_0}^{v_{k-1}} \cdots \int_{x_0}^{v_1} \frac{\hat{g}_n(v) - \sum_{j=0}^{k-1}(v-x_0)^j/j! g_0^{(j)}(x_0)}{\hat{g}_n(v)}$$
$$dv\, dv_1 \cdots dv_{k-1}$$
$$+ \widehat{A}_{k-1,n} t^{k-1} + \cdots + \widehat{A}_{0,n},$$

where

$$\widehat{A}_{j,n} = -\frac{n^{(2k-j)r_k}}{(k-1)! j!} g_0(x_0) \bigg( \widehat{H}_n^{(j)}(x_0) - \frac{(k-1)!}{(k-j)!} x_0^{k-j} \bigg), \qquad j = 0, \ldots, k-1.$$

In the following lemma, we will give the asymptotic distribution of the local processes $\mathbb{Y}_n^{\text{loc}}$ and $\widehat{\mathbb{Y}}_n^{\text{loc}}$ in terms of the $(k-1)$-fold integral of two-sided Brownian motion, $g_0(x_0)$, and $g_0^{(k)}(x_0)$ assuming that the true density $g_0$ is $k$-times continuously differentiable at $x_0$. We denote by $\bar{\mathbb{Y}}_n^{\text{loc}}$ either $\mathbb{Y}_n^{\text{loc}}$ or $\widehat{\mathbb{Y}}_n^{\text{loc}}$.

LEMMA 3.2. *Let $x_0$ be a point where $g_0$ is continuously $k$-times differentiable in a neighborhood of $x_0$ with $(-1)^k g_0^{(k)}(x_0) > 0$. Then $\bar{\mathbb{Y}}_n^{\text{loc}} \Rightarrow Y_{a,\sigma}$ in $C[-K,K]$ for each $K > 0$ where*

$$Y_{a,\sigma}(t) = \begin{cases} \sigma \int_0^t \int_0^{s_{k-1}} \cdots \int_0^{s_2} W(s_1)\, ds_1 \cdots ds_{k-1} + a(-1)^k \frac{k!}{(2k)!} t^{2k}, & t \geq 0, \\ \sigma \int_t^0 \int_{s_{k-1}}^0 \cdots \int_{s_2}^0 W(s_1)\, ds_1 \cdots ds_{k-1} + a(-1)^k \frac{k!}{(2k)!} t^{2k}, & t < 0, \end{cases}$$

where $W$ is standard two-sided Brownian motion starting at 0, $\sigma = \sqrt{g_0(x_0)}$ and $a = (-1)^k g_0^{(k)}(x_0)/k!$.

Now, let $\bar{H}_n^{\text{loc}}$ denote either $\tilde{H}_n^{\text{loc}}$ or $\widehat{H}_n^{\text{loc}}$.

LEMMA 3.3. *The localized processes $\bar{\mathbb{Y}}_n^{\text{loc}}$ and $\bar{H}_n^{\text{loc}}$ satisfy*

$$\bar{H}_n^{\text{loc}}(t) - \bar{\mathbb{Y}}_n^{\text{loc}}(t) \geq 0 \qquad \text{for all } t \geq 0,$$

*with equality if $x_0 + tn^{-1/(2k+1)}$ is a jump point of $\bar{g}_n^{(k-1)}$.*

LEMMA 3.4. *The limit process $Y_{a,\sigma}$ in Lemma 3.2 satisfies*

$$Y_{a,\sigma}(t) \stackrel{d}{=} \frac{1}{s_1} Y_{1,1}\left(\frac{t}{s_2}\right),$$

*where*

(3.2) $$s_1 = \frac{1}{\sqrt{g_0(x_0)}} \left(\frac{(-1)^k g_0^{(k)}(x_0)}{k!\sqrt{g_0(x_0)}}\right)^{(2k-1)/(2k+1)},$$

(3.3) $$s_2 = \left(\frac{\sqrt{g_0(x_0)}}{(-1)^k g_0^{(k)}(x_0)/k!}\right)^{2/(2k+1)}.$$

To show that the derivatives of $\bar{H}_n^{\text{loc}}$ are tight, we need the following lemma.

LEMMA 3.5. *For all $j \in \{0, \ldots, k-1\}$, let $\bar{A}_{j,n}$ denote either $\tilde{A}_{j,n}$ or $\widehat{A}_{j,n}$. If conjecture (1.7) holds, then*

(3.4) $$\bar{A}_{j,n} = O_p(1).$$

We now rescale the processes $\bar{\mathbb{Y}}_n^{\text{loc}}$ and $\bar{H}_n^{\text{loc}}$ so that the rescaled $\bar{\mathbb{Y}}_n^{\text{loc}}$ converges to the canonical limit process $Y_k$ defined in Lemma 3.4. Since the scaling of $\bar{\mathbb{Y}}_n^{\text{loc}}$ will be exactly the same as the one we used for $Y_k$, we define $\bar{H}_n^l$ and $\bar{\mathbb{Y}}_n^l$ by

$$\bar{H}_n^l(t) = s_1 \bar{H}_n^{\text{loc}}(s_2 t), \qquad \bar{\mathbb{Y}}_n^l(t) = s_1 \bar{\mathbb{Y}}_n^{\text{loc}}(s_2 t),$$

where $s_1$ and $s_2$ are given by (3.2) and (3.3), respectively.

LEMMA 3.6. *Let $c > 0$. Then*

$$((\bar{H}_n^l)^{(0)}, (\bar{H}_n^l)^{(1)}, \ldots, (\bar{H}_n^l)^{(2k-1)}) \Rightarrow (H_k^{(0)}, H_k^{(1)}, \ldots, H_k^{(2k-1)})$$

*in $(D[-c,c])^{2k}$, where $H_k$ is the stochastic process defined in Theorem 3.1.*

To keep this paper to a reasonable length, proofs of the results of Section 3.2 and of the main convergence Theorem 3.2 can be found in Balabdaoui and Wellner [4], Appendix 1. The arguments there are constructed along the lines of Groeneboom, Jongbloed and Wellner [18]. However, those arguments had to be adapted and further developed to be able to treat $k$-monotonicity for an arbitrary integer $k \geq 2$. In this general case, we found that it is very useful to consider perturbation functions to learn about the asymptotic behavior of the estimators. Such perturbation functions need, of course, to be permissible, that is, the resulting perturbed function must belong to the $k$-monotone class, but they also need to have a compact support to suit the local nature of the current estimation problem. It turns out that choices are rather limited and that B-splines with degree $k-1$ and support $[\tau_{n,1}, \tau_{n,k+1}]$, where $\tau_{n,1}, \ldots, \tau_{n,k+1}$ are knots of either the LSE or MLE in the neighborhood of $x_0$, are found to be the most sensible perturbation functions to consider. For a definition of B-splines, see, for example, Nürnberger [34], Theorem 2.2. For technical details on the use of B-splines for constructing perturbations, see, for example, Proposition 6.1 in Balabdaoui and Wellner [4], Appendix 1.

**4. The gap problem—spline connection.** Recall that it was assumed that $g_0$ is $k$-times continuously differentiable at $x_0$ and that $(-1)^k g_0^{(k)}(x_0) > 0$. Under a weaker assumption, Balabdaoui and Wellner [2] proved strong consistency of the $(k-1)$st derivative of the MLE and LSE. This consistency result and the above assumptions collectively imply that the number of jump points of this derivative, in a small neighborhood of $x_0$, diverges to infinity almost surely as the sample size $n \to \infty$. This "clustering" phenomenon is one of the most crucial elements in studying the local asymptotics of the estimators. The jump points then form a sequence that converges to $x_0$ almost surely and therefore the distance between two successive jump points, for example, located just before and after $x_0$, converges to 0 as $n \to \infty$. But it is not enough to know that the "gap" between these points converges to 0: an upper bound for this rate of convergence is needed.

To prove Lemma 3.1, we will focus first on the LSE because it is somewhat easier to handle through the simple form of its characterization. The arguments for the MLE could be built upon those used for the LSE, but in this case one has to deal with some extra difficulties due to the nonlinear nature of its characterization.

We start by describing the difficulties of establishing this result for the general case $k > 2$.

4.1. *Fundamental differences.* Let $\tau_n^-$ and $\tau_n^+$ be the last and first jump points of the $(k-1)$st derivative of the LSE $\tilde{g}_n$, located before and after $x_0$, respectively. To obtain a better understanding of the gap problem, we describe the reasoning used by Groeneboom, Jongbloed and Wellner [18] in order to prove that $\tau_n^+ - \tau_n^- = O_p(n^{-1/5})$ for the special case $k = 2$. The characterization of the estimator is given by

$$\text{(4.1)} \qquad \tilde{H}_n(x) \begin{cases} \geq \mathbb{Y}_n(x), & x \geq 0, \\ = \mathbb{Y}_n(x), & \text{if } x \text{ is a jump point of } \tilde{g}'_n, \end{cases}$$

where $\tilde{H}_n(x) = \int_0^x (x-t)\tilde{g}_n(t)\,dt$ and $\mathbb{Y}_n(x) = \int_0^x \mathbb{G}_n(t)\,dt$. On the interval $[\tau_n^-, \tau_n^+)$, the function $\tilde{g}'_n$ is constant since there are no more jump points in this interval. This implies that $\tilde{H}_n$ is a polynomial of degree 3 on $[\tau_n^-, \tau_n^+)$. But from the characterization in (4.1), it follows that

$$\tilde{H}_n(\tau_n^\pm) = \mathbb{Y}_n(\tau_n^\pm), \qquad \tilde{H}'_n(\tau_n^\pm) = \mathbb{Y}'_n(\tau_n^\pm).$$

These four boundary conditions allow us to fully determine the cubic polynomial $\tilde{H}_n$ on $[\tau_n^-, \tau_n^+]$. Using the explicit expression for $\tilde{H}_n$ and evaluating it at the midpoint $\bar{\tau} = (\tau_n^- + \tau_n^+)/2$, Groeneboom, Jongbloed and Wellner [18] established that

$$\tilde{H}_n(\bar{\tau}_n) = \frac{\mathbb{Y}_n(\tau_n^-) + \mathbb{Y}_n(\tau_n^+)}{2} - \frac{\mathbb{G}_n(\tau_n^+) - \mathbb{G}_n(\tau_n^-)}{8}(\tau_n^+ - \tau_n^-).$$

Groeneboom, Jongbloed and Wellner [18] refer to this as the "midpoint property." By applying the first condition (the inequality condition) in (4.1), it follows that

$$\frac{\mathbb{Y}_n(\tau_n^-) + \mathbb{Y}_n(\tau_n^+)}{2} - \frac{\mathbb{G}_n(\tau_n^+) - \mathbb{G}_n(\tau_n^-)}{8}(\tau_n^+ - \tau_n^-) \geq \mathbb{Y}_n(\bar{\tau}_n).$$

The inequality in the last display can be rewritten as

$$\frac{Y_0(\tau_n^-) + Y_0(\tau_n^+)}{2} - \frac{G_0(\tau_n^+) - G_0(\tau_n^-)}{8}(\tau_n^+ - \tau_n^-) \geq \mathbb{E}_n,$$

where $G_0$ and $Y_0$ are the true counterparts of $\mathbb{G}_n$ and $\mathbb{Y}_n$, respectively, and $\mathbb{E}_n$ is a random error. Using empirical process theory, Groeneboom, Jongbloed and Wellner [18] showed that

$$\text{(4.2)} \qquad |\mathbb{E}_n| = O_p(n^{-4/5}) + o_p((\tau_n^+ - \tau_n^-)^4).$$

On the other hand, Groeneboom, Jongbloed and Wellner [18] established that there exists a universal constant $C > 0$ such that

$$\text{(4.3)} \qquad \begin{aligned} &\frac{Y_0(\tau_n^-) + Y_0(\tau_n^+)}{2} - \frac{G_0(\tau_n^+) - G_0(\tau_n^-)}{8}(\tau_n^+ - \tau_n^-) \\ &\quad = -Cg_0''(x_0)(\tau_n^+ - \tau_n^-)^4 + o_p((\tau_n^+ - \tau_n^-)^4). \end{aligned}$$

Combining the results in (4.2) and (4.3), it follows that

$$\tau_n^+ - \tau_n^- = O_p(n^{-1/5}).$$

The problem has two main features that make the above arguments work. First, the polynomial $\tilde{H}_n$ can be fully determined on $[\tau_n^-, \tau_n^+]$ and can therefore be evaluated at any point between $\tau_n^-$ and $\tau_n^+$. Second, it can be expressed via the empirical process $\mathbb{Y}_n$ and that enables us to "get rid of" terms depending on $\tilde{g}_n$ whose rate of convergence is still unknown at this stage. We should also add that the problem is symmetric about $\bar{\tau}_n$, a property that helps in establishing the formula derived in (4.3).

When $k > 2$, it follows from the characterization of the LSE given in (2.4) that for any two successive jump points of $\tilde{g}_n^{(k-1)}$, $\tau_n^-$, $\tau_n^+$, the four equalities

$$\tilde{H}_n(\tau_n^\pm) = \mathbb{Y}_n(\tau_n^\pm) \quad \text{and} \quad \tilde{H}_n'(\tau_n^\pm) = \mathbb{Y}_n'(\tau_n^\pm)$$

still hold. However, these equations are not enough to determine the polynomial $\tilde{H}_n$, now of degree $2k - 1$, on the interval $[\tau_n^-, \tau_n^+]$. One would need $2k$ conditions to be able to achieve this. [We would be in this situation if we had equality of the higher derivatives of $\tilde{H}_n$ and $\mathbb{Y}_n$ at $\tau_n^-$ and $\tau_n^+$, i.e.,

$$(4.4) \qquad \tilde{H}_n^{(j)}(\tau_n^-) = \mathbb{Y}_n^{(j)}(\tau_n^-), \qquad \tilde{H}_n^{(j)}(\tau_n^+) = \mathbb{Y}_n^{(j)}(\tau_n^+),$$

for $j = 0, \ldots, k-1$, but the characterization (2.4) does not give this much.] Thus it becomes clear that two jump points are not sufficient to determine the piecewise polynomial $\tilde{H}_n$. However, if we consider $p > 2$ jump points $\tau_{n,0} < \cdots < \tau_{n,p-1}$ (all located, e.g., after $x_0$), then $\tilde{H}_n$ is a spline of degree $2k - 1$ with interior knots $\tau_{n,1}, \ldots, \tau_{n,p-2}$, that is, $\tilde{H}_n$ is a polynomial of degree $2k - 1$ on $(\tau_{n,j}, \tau_{n,j+1})$ for $j = 0, \ldots, p-2$ and is $(2k-2)$-times differentiable at its knot points $\tau_{n,0}, \ldots, \tau_{n,p-1}$. In the next subsection, we prove that if $p = 2k - 2$, the spline $\tilde{H}_n$ is completely determined on $[\tau_{n,0}, \tau_{n,2k-3}]$ by the conditions

$$(4.5) \qquad \tilde{H}_n(\tau_{n,i}) = \mathbb{Y}_n(\tau_{n,i}) \quad \text{and} \quad \tilde{H}_n'(\tau_{n,i}) = \mathbb{Y}_n'(\tau_{n,i}),$$

$$i = 0, \ldots, 2k-3.$$

This result proves to be very useful for determining the stochastic order of the distance between two successive jump points in a small neighborhood of $x_0$ if conjecture (1.7) on the uniform boundedness of the error in the "nonclassical" Hermite interpolation problem via splines of odd degree defined in (1.6) holds.

4.2. *The gap problem for the LSE—Hermite interpolation.* In the next lemma, we prove that given $2k - 2$ successive jump points $\tau_{n,0} < \cdots < \tau_{n,2k-3}$ of $\tilde{g}_n^{(k-1)}$, $\tilde{H}_n$ is the unique solution of the Hermite problem given by (4.5). In the following, we will omit writing the subscript $n$ explicitly in the knots, but their dependence on the sample size should be kept in mind.

LEMMA 4.1. *The function $\tilde{H}_n$ characterized by (2.4) is a spline of degree $2k-1$. Moreover, given any $2k-2$ successive jump points of $\tilde{H}_n^{(2k-1)}$, $\tau_0 < \cdots < \tau_{2k-3}$, the $(2k-1)$st spline $\tilde{H}_n$ is uniquely determined on $[\tau_0, \tau_{2k-3}]$ by the values of the process $\mathbb{Y}_n$ and of its derivative $\mathbb{Y}'_n$ at $\tau_0, \ldots, \tau_{2k-3}$.*

PROOF. We know that for any jump point $\tau$ of $\tilde{H}_n^{(2k-1)}$, we have

$$\tilde{H}_n(\tau) = \mathbb{Y}_n(\tau) \quad \text{and} \quad \tilde{H}'_n(\tau) = \mathbb{Y}'_n(\tau).$$

This can be viewed as a *Hermite interpolation problem* if we consider that the *interpolated function* is the process $\mathbb{Y}_n$ and that the *interpolating spline* is $\tilde{H}_n$ (see, e.g., Nürnberger [34], Definition 3.6, pages 108 and 109). Existence and uniqueness of the spline interpolant follows easily from the Schoenberg–Whitney–Karlin–Ziegler theorem (Schoenberg and Whitney [37], Theorem 3, page 258; Karlin and Ziegler [23], Theorem 3, page 529; Nürnberger [34], Theorem 3.7, page 109; DeVore and Lorentz [9], Theorem 9.2, page 162). □

In the following lemma, we prove a preparatory result that will be used later for deriving the stochastic order of the distance between successive knots, $\tau_0, \ldots, \tau_{2k-3}$, of $\tilde{g}_n$ in a neighborhood of $x_0$. Let $\mathcal{H}_k$ again denote the spline interpolation operator which assigns to each differentiable function $f$ the unique spline $\mathcal{H}_k[f]$ with interior knots $\tau_1, \ldots, \tau_{2k-4}$ and degree $2k-1$, and satisfying the boundary conditions given in (1.6).

LEMMA 4.2. *Let $\bar{\tau} \in \bigcup_{i=0}^{2k-4} (\tau_i, \tau_{i+1})$. If $e_k(t)$ denotes the error at $t$ of the Hermite interpolation of the function $x^{2k}/(2k)!$, that is,*

$$e_k(t) = \frac{t^{2k}}{(2k)!} - \mathcal{H}_k\!\left[\frac{x^{2k}}{(2k)!}\right]\!(t),$$

*then*

(4.6) $$g_0^{(k)}(\bar{\tau}) e_k(\bar{\tau}) \leq \mathbb{E}_n + \mathbb{R}_n,$$

*where $\mathbb{E}_n$ defined in (4.8) is a random error and $\mathbb{R}_n$ defined in (4.9) is a remainder that both depend on the knots $\tau_0, \ldots, \tau_{2k-3}$ and the point $\bar{\tau}$.*

PROOF. Let $\bar{\tau} \in \bigcup_{i=0}^{2k-4}(\tau_i, \tau_{i+1})$. From the characterization in (2.4) and the fact that $\tilde{H}_n = \mathcal{H}_k[\mathbb{Y}_n]$ on $[\tau_0, \tau_{2k-3}]$, it follows that

$$\mathcal{H}_k[\mathbb{Y}_n](\bar{\tau}) \geq \mathbb{Y}_n(\bar{\tau}).$$

Let $Y_0$ be the true counterpart of $\mathbb{Y}_n$, that is, $Y_0(x) = \int_0^x (x-t)^{k-1} g_0(t)\,dt/(k-1)!$. We can then rewrite the previous inequality as

(4.7) $$\mathcal{H}_k[Y_0](\bar{\tau}) - Y_0(\bar{\tau}) \geq -\mathbb{E}_n(\bar{\tau}),$$

where

(4.8) $$\mathbb{E}_n = \mathcal{H}_k[\mathbb{Y}_n - Y_0](\bar{\tau}) - [\mathbb{Y}_n - Y_0](\bar{\tau}).$$

Based on the working assumptions, the function $Y_0$ is $(2k)$-times continuously differentiable in a small neighborhood of $x_0$. Now, Taylor expansion of $Y_0(t)$ with integral remainder around $\bar{\tau}$ up to the order $2k$ yields

$$Y_0(t) = \sum_{j=0}^{2k-1} \frac{(t-\bar{\tau})^j}{j!} Y_0^{(j)}(\bar{\tau}) + \int_{\bar{\tau}}^{\tau_{2k-3}} \frac{(t-u)_+^{2k-1}}{(2k-1)!} g_0^{(k)}(u)\, du,$$

for all $t \in [\tau_0, \tau_{2k-3}]$. Using this expansion, along with the fact that the operator $\mathcal{H}_k$ is linear and preserves polynomials of degree $2k - 1$, we can rewrite the inequality in (4.7) as

$$\frac{1}{(2k-1)!} \int_{\bar{\tau}}^{\tau_{2k-3}} \mathcal{H}_k[(t-u)_+^{2k-1}](\bar{\tau}) g_0^{(k)}(u)\, du \geq -\mathbb{E}_n.$$

In the previous display, $\mathcal{H}_k[(t-u)_+^{2k-1}](\bar{\tau})$ is the Hermite spline interpolant of the truncated power function $t \mapsto (t-u)_+^{2k-1}$ ($u$ is fixed), evaluated at the point $\bar{\tau}$. Now, we can rewrite the left-hand side of the previous inequality as

$$\int_{\bar{\tau}}^{\tau_{2k-3}} \frac{1}{(2k-1)!} \mathcal{H}_k[(t-u)_+^{2k-1}](\bar{\tau}) g_0^{(k)}(u)\, du$$

$$= g_0^{(k)}(\bar{\tau}) \frac{1}{(2k-1)!} \int_{\bar{\tau}}^{\tau_{2k-3}} \mathcal{H}_k[(t-u)_+^{2k-1}](\bar{\tau})\, du$$

(4.9)

$$+ \frac{1}{(2k-1)!} \int_{\bar{\tau}}^{\tau_{2k-3}} \mathcal{H}_k[(t-u)_+^{2k-1}](\bar{\tau})(g_0^{(k)}(u) - g_0^{(k)}(\bar{\tau}))\, du$$

$$= g_0^{(k)}(\bar{\tau}) \frac{1}{(2k-1)!} \mathcal{H}_k\left[\int_{\bar{\tau}}^{\tau_{2k-3}} [(t-u)_+^{2k-1}]\, du\right](\bar{\tau}) + \mathbb{R}_n,$$

once again using linearity of the operator $\mathcal{H}_k$. The remainder $\mathbb{R}_n$ is equal to the Hermite interpolant of the function

$$t \mapsto \frac{1}{(2k-1)!} \int_{\bar{\tau}}^{t} \frac{(t-u)^{2k-1}}{(2k-1)!} (g_0^{(k)}(u) - g_0^{(k)}(\bar{\tau}))\, du$$

at the point $\bar{\tau}$. On the other hand, we can further rewrite the integral term in (4.9) as

$$\frac{1}{(2k-1)!} \mathcal{H}_k\left[\int_{\bar{\tau}}^{\tau_{2k-3}} (t-u)_+^{2k-1}\, du\right](\bar{\tau})$$

$$= \frac{1}{(2k-1)!} \mathcal{H}_k\left[\int_{\bar{\tau}}^{t} (t-u)^{2k-1}\, du\right](\bar{\tau})$$

$$= \frac{1}{(2k)!} \mathcal{H}_k[(t-\bar{\tau})^{2k}](\bar{\tau}).$$

In other words, the integral term in (4.9) is nothing but the value of the Hermite spline interpolant of the function $t \mapsto (t-\bar{\tau})^{2k}/(2k)!$ at the point $\bar{\tau}$. As claimed in the lemma, this value is also equal to $-e_k(\bar{\tau})$, where $e_k$ is the error of the Hermite interpolation of the function $x^{2k}/(2k)!$. Indeed, let $P_{2k-1}(t) = (t-\bar{\tau})^{2k}/(2k)! - t^{2k}/(2k)!$. Since $P_{2k-1}$ is a polynomial of degree $2k-1$, we have

$$\mathcal{H}_k\left[\frac{(x-\bar{\tau})^{2k}}{(2k)!}\right](t) = \mathcal{H}_k\left[\frac{x^{2k}}{(2k)!}\right](t) + P_{2k-1}(t).$$

If $t = \bar{\tau}$, $P_{2k-1}(\bar{\tau}) = 0 - \bar{\tau}^{2k}/(2k)! = -\bar{\tau}^{2k}/(2k)!$, which implies that

$$\mathcal{H}_k\left[\frac{(x-\bar{\tau})^{2k}}{(2k)!}\right](\bar{\tau}) = \mathcal{H}_k\left[\frac{x^{2k}}{(2k)!}\right](\bar{\tau}) - \frac{\bar{\tau}^{2k}}{(2k)!} = -e_k(\bar{\tau}). \qquad \Box$$

The error $e_k$ defined in Lemma 4.2 can be recognized as a monospline of degree $2k$ with $2k-2$ simple knots $\tau_0, \ldots, \tau_{2k-3}$. For a definition of monosplines, see, for example, Micchelli [32], Bojanov, Hakopian and Sahakian [5], Nürnberger [34], page 194, or DeVore and Lorentz [9], page 136. In the next lemma, we state an important property of $e_k$.

LEMMA 4.3. *The function $x \mapsto e_k(x)$ has no zeros other than $\tau_0, \ldots, \tau_{2k-3}$ in $[\tau_0, \tau_{2k-3}]$. Furthermore, $(-1)^k e_k \geq 0$ on $[\tau_0, \tau_{2k-3}]$.*

PROOF. See Balabdaoui and Wellner [4], Appendix 3. $\Box$

In Lemma 4.2, the key inequality in (4.6) can be rewritten as

(4.10) $$(-1)^k g_0^{(k)}(\bar{\tau}) \cdot (-1)^k e_k(\bar{\tau}) \leq \mathbb{E}_n + \mathbb{R}_n,$$

where the first factor on the right-hand side is already known to be positive by $k$-monotonicity of $g_0$. Lemmas 4.4 and 4.5 are the final steps toward establishing the order of the gap for the LSE based on conjecture (1.7).

LEMMA 4.4. *If conjecture (1.7) holds, then $\mathbb{E}_n$ in (4.6) of Lemma 4.2 satisfies*

$$|\mathbb{E}_n| = O_p(n^{-k/(2k+1)}) + o_p((\tau_{2k-3} - \tau_0)^{2k}).$$

PROOF. We have

$$\mathbb{E}_n = \mathcal{H}_k[\mathbb{Y}_n - Y_0](\bar{\tau}) - [\mathbb{Y}_n - Y_0](\bar{\tau}).$$

Using (generalized) Taylor expansions of $\mathbb{Y}_n$ and $Y_0$ around the point $\bar{\tau}$ up to order $k-1$ yields

$$\mathbb{Y}_n(t) - Y_0(t) = \sum_{j=0}^{k-1} \frac{(t-\bar{\tau})^j}{j!}[\mathbb{Y}_n^{(j)}(\bar{\tau}) - Y_0^{(j)}(\bar{\tau})] + \int_{\bar{\tau}}^{t} \frac{(t-x)^{k-1}}{(k-1)!} \, d(\mathbb{G}_n - G_0)(x),$$

and therefore,

$$\begin{aligned}
\mathbb{E}_n &= \mathcal{H}_k\bigg[\int_{\bar{\tau}}^t \frac{1}{(k-1)!}(t-x)^{k-1}\,d(\mathbb{G}_n - G_0)(x)\bigg](\bar{\tau}) \\
&= \mathcal{H}_k\bigg[\int_{\bar{\tau}}^{\tau_{2k-3}} g_t(x)\,d(\mathbb{G}_n - G_0)(x)\bigg](\bar{\tau}), \qquad \text{where } g_t(x) = \frac{(t-x)_+^{k-1}}{(k-1)!} \\
&= \int_{\bar{\tau}}^{\tau_{2k-3}} \mathcal{H}_k[g_t(x)](\bar{\tau})\,d(\mathbb{G}_n - G_0)(x), \qquad \text{by linearity of } \mathcal{H}_k \\
&= \int_{\tau_0}^{\tau_{2k-3}} f_{\bar{\tau}}(x)\,d(\mathbb{G}_n - G_0)(x).
\end{aligned}$$

Given $x \in [\bar{\tau}, \tau_{2k-3}]$, $f_{\bar{\tau}}(x) = \mathcal{H}_k[g_t(x)](\bar{\tau})1_{[\bar{\tau},\tau_{2k-3}]}(x)$, where $\mathcal{H}_k[g_t(x)](\bar{\tau})$ is the value at $\bar{\tau}$ of the Hermite spline interpolant of the function $t \mapsto g_t(x) = (t-x)_+^{k-1}/(k-1)!$. Thus $f_{\bar{\tau}}(x)$ depends on the knots $\tau_0, \ldots, \tau_{2k-3}$ and the point $s = \bar{\tau} \in [\tau_0, \tau_{2k-3}]$ and can be viewed as an element of the class of functions

$$\begin{aligned}
(4.11) \quad \mathcal{F}_{y_0,R} = \{f_s(x) &= f_{s,y_0,\ldots,y_{2k-3}}(x) : x \in [y_0, y_{2k-3}], s \in [y_0, y_{2k-3}], \\
x_0 - \delta &\leq y_0 < y_1 < \cdots < y_{2k-3} \leq y_0 + R \leq x_0 + \delta\},
\end{aligned}$$

where $\delta > 0$ is a fixed small number. In view of conjecture (1.7), together with the triangle inequality, there exists a constant $C > 0$ depending only on $k$ such that

$$|f_s(x)| \leq C(y_{2k-3} - y_0)^{k-1}1_{[y_0, y_{2k-3}]}(x)$$

and hence the collection $\mathcal{F}_{y_0,R}$ has envelope function $F_{y_0,R}$ given by

$$F_{y_0,R}(x) = CR^{k-1}1_{[y_0, y_0+R]}(x).$$

Furthermore, $\mathcal{F}_{y_0,R}$ is a VC-subgraph collection of functions (see Proposition A.1 in the Appendix for a detailed argument) and hence by van der Vaart and Wellner [40], Theorem 2.6.7, page 141,

$$\sup_Q N(\varepsilon\|F\|_{Q,2}, \mathcal{F}_{y_0,R}, L_2(Q)) \leq \left(\frac{K}{\varepsilon}\right)^{V_k},$$

for $0 < \varepsilon < 1$, where $V_k = 2(V(\mathcal{F}_{y_0,R}) - 1)$ with $V(\mathcal{F}_{y_0,R})$ the VC-dimension of the collection of subgraphs and where the constant $K$ depends only on $V(\mathcal{F}_{y_0,R})$ [note that from our proof of Proposition A.1, it is clear that $V(\mathcal{F}_{y_0,R})$ depends only of $k$]. It follows that

$$\sup_Q \int_0^1 \sqrt{1 + \log N(\varepsilon\|F_{y_0,R}\|_{Q,2}, \mathcal{F}_{y_0,R}, L_2(Q))}\,d\varepsilon < \infty.$$

On the other hand,

$$EF_{y_0,R}^2(X_1) = C^2 R^{2(k-1)} \int_{y_0}^{y_0+R} g_0(x)\,dx \leq C^2 M R^{2k-1}$$

with $M \equiv g_0(x_0 - \delta)$. Application of Lemma A.1 with $d = k$ yields

$$|\mathbb{E}_n| = o_p((\tau_{2k-3} - \tau_0)^{2k}) + O_p(n^{-2k/(2k+1)}). \qquad \square$$

LEMMA 4.5. *If the bound in (1.8) holds, then $\mathbb{R}_n$ of Lemma 4.2 satisfies*

$$|\mathbb{R}_n| = o_p((\tau_{2k-3} - \tau_0)^{2k}).$$

PROOF. By definition, $\mathbb{R}_n$ is the value at $\bar{\tau}$ of the Hermite spline interpolant of the function

$$(4.12) \qquad t \mapsto \int_{\bar{\tau}}^{t} \frac{(t-u)^{2k-1}}{(2k-1)!}(g_0^{(k)}(u) - g_0^{(k)}(\bar{\tau}))\,du.$$

By (1.8), there exists a constant $D > 0$ depending only on $k$ such that

$$|\mathbb{R}_n| \leq D \sup_{t \in [\tau_0, \tau_{2k-3}]} |g_0^{(k)}(t) - g_0^{(k)}(\bar{\tau})|(\tau_{2k-3} - \tau_0)^{2k}.$$

In the previous bound, we used the fact that the $(2k)$-times derivative of the function in (4.12) is $g_0^{(k)}(t) - g_0^{(k)}(\bar{\tau})$. But, note that this derivative is $o_p(1)$, which follows from uniform continuity of $g_0^{(k)}$ on compacts. This, in turn, implies the claimed bound. $\square$

PROOF OF LEMMA 3.1 FOR THE LSE. Let $j_0 \in \{0,\ldots,2k-4\}$ be such that $[\tau_{j_0}, \tau_{j_0+1}]$ is the largest knot interval, that is, $\tau_{j_0+1} - \tau_{j_0} = \max_{0 \leq j \leq 2k-4}(\tau_{j+1} - \tau_j)$. Let $a = \tau_0$, $b = \tau_{2k-3}$. Using the inequality in (4.10) and noting that the bounds on $\mathbb{R}_n$ and $\mathbb{E}_n$ are independent of the choice of $\bar{\tau}$ in $\bigcup_{j=0}^{2k-4}(\tau_j, \tau_{j+1})$, it follows that

$$\sup_{\bar{\tau} \in (\tau_{j_0}, \tau_{j_0+1})} (-1)^k e_k(\bar{\tau}) \leq O_p(n^{-2k/(2k+1)}) + o_p((\tau_{2k-3} - \tau_0)^{2k}).$$

Now, on the interval $[\tau_{j_0}, \tau_{j_0+1}]$, the Hermite spline interpolant of the function $x^{2k}/(2k)!$ reduces to a polynomial of degree $2k - 1$. On the other hand, the best uniform approximation of the function $x^{2k}$ on $[\tau_{j_0}, \tau_{j_0+1}]$ from the space of polynomials of degree $\leq 2k - 1$ is given by the polynomial

$$(4.13) \quad x \mapsto x^{2k} - \left(\frac{\tau_{j_0+1} - \tau_{j_0}}{2}\right)^{2k} \frac{1}{2^{2k-1}} T_{2k}\left(\frac{2x - (\tau_{j_0} + \tau_{j_0+1})}{\tau_{j_0+1} - \tau_{j_0}}\right),$$

where $T_{2k}$ is the Chebyshev polynomial of degree $2k$ (defined on $[-1,1]$); see, for example, Nürnberger [34], Theorem 3.23, page 46, or DeVore and Lorentz [9], Theorem 6.1, page 75. It follows that

$$(4.14) \quad \sup_{\bar{\tau} \in (\tau_{j_0}, \tau_{j_0+1})} (-1)^k e_k(\bar{\tau}) \geq \left\| \frac{T_{2k}}{2^{4k-1}(2k)!} \right\|_\infty (\tau_{j_0+1} - \tau_{j_0})^{2k}$$

$$= \frac{1}{2^{4k-1}(2k)!} (\tau_{j_0+1} - \tau_{j_0})^{2k}$$

since $\|T_{2k}\|_\infty = 1$. But

$$\tau_{2k-3} - \tau_0 = \sum_{j=0}^{2k-4} (\tau_{j+1} - \tau_j) \leq (2k-3)(\tau_{j_0+1} - \tau_{j_0}).$$

Hence,

$$\sup_{\bar{\tau} \in (\tau_{j_0}, \tau_{j_0+1})} (-1)^k e_k(\bar{\tau}) \geq \frac{1}{(2k-3)^{2k} 2^{4k-1}(2k)!} (\tau_{2k-3} - \tau_0)^{2k}.$$

Combining the results obtained above, we conclude that

$$\frac{(-1)^k g_0^{(k)}(x_0)}{(2k-3)^{2k} 2^{4k-1}(2k)!} (\tau_{2k-3} - \tau_0)^{2k} \leq O_p(n^{-2k/(2k+1)}) + o_p((\tau_{2k-3} - \tau_0)^{2k}),$$

which implies that $\tau_{2k-3} - \tau_0 = O_p(n^{-1/(2k+1)})$. □

4.3. *The gap problem for the MLE.* To show Lemma 3.1 for the MLE, one needs to deal with an extra difficulty posed by the nonlinear form of the characterization of this estimator as given in (2.6). In the following, we show how one can get around this difficulty. The main idea is to "linearize" the characterization of the MLE and hence be able to re-use the arguments developed for the LSE in the previous subsection.

LEMMA 4.6. *Let $\tau_0, \ldots, \tau_{2k-3}$ be $2k-2$ successive jump points of $\hat{g}_n^{(k-1)}$. Then*

$$\mathcal{H}_k[\mathbb{Y}_n] - \mathbb{Y}_n \geq g_0(\tau_0)(\check{f}_n - \mathcal{H}_k[\check{f}_n] + \Delta_n - \mathcal{H}_k[\Delta_n])$$

*on $[\tau_0, \tau_{2k-3}]$, where $\mathbb{Y}_n$ is the same empirical process introduced in (2.2),*

$$\check{f}_n(x) \equiv -\int_{\tau_0}^t \frac{(x-t)^{k-1}}{(k-1)!} \left( \frac{1}{\hat{g}_n(t)} - \frac{1}{g_0(\tau_0)} \right) d(\hat{G}_n(t) - G_0(t))$$

*and*

$$\Delta_n(x) \equiv \int_{\tau_0}^x \frac{(x-t)^{k-1}}{(k-1)!} \left( \frac{1}{\hat{g}_n(t)} - \frac{1}{g_0(\tau_0)} \right) d(\mathbb{G}_n(t) - G_0(t)).$$

PROOF. Let $\hat{G}_n(x) = \int_0^x \hat{g}_n(s)\,ds$. The characterization in (2.6) can be rewritten as

$$(4.15) \quad \int_0^x \frac{(x-t)^{k-1}}{\hat{g}_n(t)}\,d(\hat{G}_n(t) - \mathbb{G}_n(t)) \begin{cases} \geq 0, & \text{for } x > 0, \\ = 0, & \text{if } x \text{ is a jump point of } \hat{g}_n^{(k-1)}. \end{cases}$$

Note that when $x$ is a jump point of $\hat{g}_n^{(k-1)}$, the two parts of (4.15) imply that the first derivative of the function on the right-hand side is equal to 0 at the jump point $x$, that is,

$$(4.16) \quad \int_0^x \frac{(x-t)^{k-2}}{\hat{g}_n(t)}\,d(\hat{G}_n(t) - \mathbb{G}_n(t)) = 0.$$

For $x > 0$, let

$$\hat{H}_n(x) = \int_0^x \frac{(x-t)^{k-1}}{(k-1)!}\,d\hat{G}_n(t).$$

Note that $\hat{H}_n \neq \widetilde{H}_n$ defined in (2.5) and that on $[\tau_0, \tau_{2k-3}]$, $\hat{H}_n$ is a spline of degree $2k-1$ with knots $\tau_0, \ldots, \tau_{2k-3}$. For $x \in [\tau_0, \tau_{2k-3}]$, we can write

$$\int_0^x \frac{(x-t)^{k-1}}{\hat{g}_n(t)}\,d(\hat{G}_n(t) - \mathbb{G}_n(t))$$
$$= \frac{1}{g_0(\tau_0)} \int_0^x (x-t)^{k-1}\,d(\hat{G}_n(t) - \mathbb{G}_n(t))$$
$$\quad + \int_0^x (x-t)^{k-1}\left(\frac{1}{\hat{g}_n(t)} - \frac{1}{g_0(\tau_0)}\right) d(\hat{G}_n(t) - \mathbb{G}_n(t))$$
$$= \frac{(\hat{H}_n(x) - \mathbb{Y}_n(x))}{g_0(\tau_0)} + \int_0^{\tau_0} (x-t)^{k-1}\left(\frac{1}{\hat{g}_n(t)} - \frac{1}{g_0(\tau_0)}\right) d(\hat{G}_n(t) - \mathbb{G}_n(t))$$
$$\quad + \int_{\tau_0}^x (x-t)^{k-1}\left(\frac{1}{\hat{g}_n(t)} - \frac{1}{g_0(\tau_0)}\right) d(\hat{G}_n(t) - G_0(t))$$
$$\quad + \int_{\tau_0}^x (x-t)^{k-1}\left(\frac{1}{\hat{g}_n(t)} - \frac{1}{g_0(\tau_0)}\right) d(G_0(t) - \mathbb{G}_n(t))$$
$$= \frac{1}{g_0(\tau_0)}(\hat{H}_n(x) - \mathbb{Y}_n(x)) + p_n(x) - \check{f}_n(x) - \Delta_n(x).$$

Note that

$$p_n(x) \equiv \int_0^{\tau_0} (x-t)^{k-1}\left(\frac{1}{\hat{g}_n(t)} - \frac{1}{g_0(\tau_0)}\right) d(\hat{G}_n(t) - \mathbb{G}_n(t))$$

is a polynomial of degree $k-1$. From (4.15) and (4.16), it follows that $\hat{H}_n$ is the Hermite spline interpolant of the function

$$\mathbb{Y}_n + g_0(\tau_0)\{-p_n + \check{f}_n + \Delta_n\}$$

such that

$$\hat{H}_n \geq \mathbb{Y}_n + g_0(\tau_0)(-p_n + \check{f}_n + \Delta_n).$$

Hence,

$$\mathcal{H}_k[\mathbb{Y}_n + g_0(\tau_0)\{-p_n + \check{f}_n + \Delta_n\}] \geq \mathbb{Y}_n + g_0(\tau_0)\{-p_n + \check{f}_n + \Delta_n\}$$

on $[\tau_0, \tau_{2k-3}]$ or, equivalently,

$$\mathcal{H}_k[\mathbb{Y}_n] - \mathbb{Y}_n \geq g_0(\tau_0)(\check{f}_n - \mathcal{H}_k[\check{f}_n] + \Delta_n - \mathcal{H}_k[\Delta_n]). \qquad \square$$

Since $\mathcal{H}_k[\mathbb{Y}_n] - \mathbb{Y}_n$ has already been studied for the purposes of proving the order of the gap in the case of the LSE, the final step is to evaluate each of the interpolation errors

(4.17) $$\mathcal{E}_1 = \check{f}_n - \mathcal{H}_k[\check{f}_n] \quad \text{and} \quad \mathcal{E}_2 = \Delta_n - \mathcal{H}_k[\Delta_n].$$

LEMMA 4.7. *Let $\mathcal{E}_1$ and $\mathcal{E}_2$ be the interpolation errors defined in (4.17). Then*

$$\|\mathcal{E}_1\|_\infty = o_p((\tau_{2k-3} - \tau_0)^{2k})$$

*and*

$$\|\mathcal{E}_2\|_\infty = o_p((\tau_{2k-3} - \tau_0)^{2k}) + O_p(n^{-2k/(2k+1)}).$$

PROOF. A detailed proof can be found in Balabdaoui and Wellner [4], Appendix 3. $\square$

PROOF OF LEMMA 3.1 FOR THE MLE. From our study of the distance between the knots of the LSE, and using very similar calculations, we can show that for all $\bar{\tau} \in \bigcup_{j=0}^{2k-4}(\tau_j, \tau_{j+1})$,

$$(-1)^k g_0^{(k)}(\bar{\tau})(-1)^k e_k(\bar{\tau}) \leq \mathbb{E}_n + \mathbb{R}_n - g_0(\tau_0)(\mathcal{E}_1(\bar{\tau}) + \mathcal{E}_2(\bar{\tau})),$$

which implies that by the results obtained for the LSE,

$$D(\tau_{2k-3} - \tau_0)^{2k}(1 + o_p(1)) \leq O_p(n^{-2k/(2k+1)}) + g_0(\tau_0)(\|\mathcal{E}_1\|_\infty + \|\mathcal{E}_2\|_\infty)$$

for some constant $D > 0$ depending on $k$ and $x_0$. Hence, it follows from Lemma 4.7 that

$$D(\tau_{2k-3} - \tau_0)^{2k}(1 + o_p(1)) \leq O_p(n^{-2k/(2k+1)}),$$

which yields the order $n^{-1/(2k+1)}$ for the distance between the knots of the MLE in the neighborhood of $x_0$. $\square$

**5. Conclusions and discussion.** As noted in Section 1, one of the motivations for this work was to try to approach the problem of pointwise limit theory for the MLE's in both the forward and inverse problems for the family of completely monotone densities on $\mathbb{R}^+$. This is one very important special case of the family of nonparametric mixture models with a smooth kernel as was mentioned in part (b) of our discussion in Section 1. Jewell [22] established consistency of the MLE's of $g \in \mathcal{D}_\infty$ and the corresponding mixing distribution function $F$ in this setting, but local rates of convergence and limiting distribution theory remain unknown. Our initial hope was that we might be able to learn about the problem with $k = \infty$ by studying the problem for fixed $k$ and then taking limits as $k \to \infty$. Unfortunately, we now believe that new tools and methods will be needed. The following discusses the state of affairs as we understand it now.

In terms of rates of convergence and localization properties, our development here shows that the local behavior of the estimators near a fixed point $x_0 > 0$ becomes dependent on an increasing number of jump points or knots in the spline problem. In other words, one needs to consider $2k - 2$ consecutive jump points (knots) $\tau_{0,n} < \cdots < \tau_{n,2k-3}$ of the $(k-1)$st derivative of the estimators in a neighborhood of $x_0$ in order to be able to find a bound on $\tau_{n,j+1} - \tau_{n,j}, j = 0, \ldots, 2k-4$, as $n \to \infty$. Thus the problem becomes increasingly "less local" with increasing $k$ and this leads us to suspect that the situation in the $k = \infty$ (or completely monotone) problem might be only "weakly local" or perhaps even "completely nonlocal" in senses yet to be precisely defined.

Another aspect of this problem is that although the MLE is asymptotically equivalent to the (mass-unconstrained) LSE for each fixed $k$ if our conjecture (1.7) holds, they seem to differ increasingly as $k$ increases. For $k = 1$, the MLE and the LSE are identical; for $k = 2$, the MLE differs from the (mass-unconstrained) LSE, but the LSE always has total mass 1. For $k \geq 3$, the MLE and LSE differ, and, moreover, the total amount of mass in the unconstrained LSE for $n = 1$ is $M_k = ((2k-1)/k)(1 - 1/(2k-1))^{k-1} \nearrow 2e^{-1/2} \approx 1.21306\ldots \neq 1$ as $k \to \infty$. We do not know how the mass of the unconstrained LSE behaves jointly in $n$ and $k$, even though (by consistency) the mass of the LSE converges to 1 as $n \to \infty$ for fixed $k$. We also do not even know if the unconstrained LSE exists for the scale mixture of exponentials, even though it is clear that the constrained estimator (defined by the least squares criterion minimized over $\mathcal{D}_k$ rather than $\mathcal{M}_k$) with mass 1 does exist. Since our current proof techniques rely so heavily on showing equivalence between the MLE and the (unconstrained) LSE, it seems likely that new methods will be required. We do not know if the (mass)-constrained LSE's and the MLE's are asymptotically equivalent either for finite $k$ or for $k = \infty$. Our current plan is to study the constrained LSE's with total mass constrained to be 1 for finite sample sizes in order to investigate the asymptotic equivalence of these mass-constrained LSE's and the MLE's and to (perhaps) extend this study to $k = \infty$ via limits on $k$. We do not yet know the "right" Gaussian version of the estimation problem in the completely monotone case.

Another way to view these difficulties might be to take the following perspective: since more knowledge is available concerning the MLE's for the families $\mathcal{D}_k$ with $k$ finite and since $\mathcal{D}_\infty$ is the intersection of all of the $\mathcal{D}_k$'s (and hence well approximated by $\mathcal{D}_k$ with $k$ large), we can fruitfully consider estimation via model selection, choosing $k$ based on the data, over the collection $\bigcup_{k=1}^\infty \mathcal{D}_k$.

In summary, we have tried to shed some more light on the local behavior of two nonparametric estimators of a $k$-monotone density, the Maximum Likelihood and Least Squares estimators. We have shown that they are both adaptive splines of degree $k-1$ with knots determined by the data and their corresponding criterion functions. When $(-1)^k g_0^{(k)}(x_0) > 0$, the distance between their knots in a neighborhood of a point $x_0 > 0$ was shown to be $n^{-1/(2k+1)}$ if a conjecture concerning the uniform boundedness of the interpolation error in a new Hermite interpolation problem holds. Once this control of the distance between the knots is available, pointwise limit distribution theory follows via a route paralleling previous results for $k = 1, 2$. Although we do not exclude the possibility that this order could be established via other approaches, we hope that the techniques developed here demonstrate that there could still be many interesting and powerful connections between statistics and approximation theory.

## APPENDIX: PROOFS FROM EMPIRICAL PROCESSES THEORY

The following proposition is a slight generalization of Lemma 4.1 of Kim and Pollard [24], page 201.

LEMMA A.1. *Let $\mathcal{F}$ be a collection of functions defined on $[x_0 - \delta, x_0 + \delta]$, with $\delta > 0$ small. Suppose that for a fixed $x \in [x_0 - \delta, x_0 + \delta]$ and $R > 0$ such that $[x, x+R] \subseteq [x_0 - \delta, x_0 + \delta]$, the collection*

$$\mathcal{F}_{x,R} = \{f_{x,y} \equiv f 1_{[x,y]}, \ f \in \mathcal{F}, \ x \leq y \leq x + R\}$$

*admits an envelope $F_{x,R}$ such that*

$$E F_{x,R}^2(X_1) \leq K R^{2d-1}, \qquad R \leq R_0,$$

*for some $d \geq 1/2$, where $K > 0$ depends only on $x_0$ and $\delta$. Moreover, suppose that*

(A.1) $$\sup_Q \int_0^1 \sqrt{\log N(\eta \|F_{x,R}\|_{Q,2}, \mathcal{F}_{x,R}, L_2(Q))} \, d\eta < \infty.$$

*Then for each $\varepsilon > 0$, there exist random variables $M_n$ of order $O_p(1)$ such that*

$$|(\mathbb{G}_n - G_0)(f_{x,y})| \leq \varepsilon |y-x|^{k+d} + n^{-(k+d)/(2k+1)} M_n \tag{A.2}$$

*for $|y - x| \leq R_0$.*

PROOF. By van der Vaart and Wellner [40], Theorem 2.14.1, page 239, it follows that

$$E\left\{\sup_{f_{x,y} \in \mathcal{F}_{x,R}} |(\mathbb{G}_n - G_0)(f_{x,R})|\right\}^2 \leq \frac{K}{n} EF_{x,R}^2(X_1) = O(n^{-1} R^{2d-1}) \tag{A.3}$$

for some constant $K > 0$ depending only on $x_0$, $\delta$ and the entropy integral in (A.1). For any $f_{x,y} \in \mathcal{F}_{x,R}$, we write

$$(\mathbb{P}_n - P_0)(f_{x,y}) = (\mathbb{G}_n - G_0)(f_{x,y})$$

and define $M_n$ by

$$M_n = \inf\{D > 0 : |(\mathbb{P}_n - P_0)(f_{x,y})| \leq \varepsilon(y-x)^{k+d} + n^{-(k+d)/(2k+1)} D,$$

$$\text{for all } f_{x,y} \in \mathcal{F}_{x,R}\}$$

and $M_n = \infty$ if no $D > 0$ satisfies the required inequality. For $1 \leq j \leq \lfloor Rn^{1/(2k+1)} \rfloor = j_n$, we have

$$P(M_n > m)$$

$$\leq P(|(\mathbb{P}_n - P_0)(f_{x,y})| > \varepsilon(y-x)^{k+d} + n^{-(k+d)/(2k+1)} m$$

$$\text{for some } f_{x,y} \in \mathcal{F}_{x,R})$$

$$\leq \sum_{1 \leq j \leq j_n} P\{n^{(k+d)/(2k+1)} |(\mathbb{P}_n - P_0)(f_{x,y})| > \varepsilon(j-1)^{k+d} + m$$

$$\text{for some } f_{x,y} \in \mathcal{F}_{x,R}, (j-1)n^{-1/(2k+1)} \leq y - x \leq jn^{-1/(2k+1)}\}$$

$$\leq \sum_{1 \leq j \leq j_n} n^{2(k+d)/(2k+1)} \frac{E\{\sup_{y:0 \leq y-x < jn^{-1/(2k+1)}} |(\mathbb{P}_n - P_0)(f_{x,y-x})|\}^2}{(\varepsilon(j-1)^{k+d} + m)^2}$$

$$= \sum_{1 \leq j \leq j_n} n^{2(k+d)/(2k+1)} \frac{E\{\sup_{f_{x,y-x} \in \mathcal{F}_{x,jn^{-1/(2k+1)}}} |(\mathbb{P}_n - P_0)(f_{x,y-x})|\}^2}{(\varepsilon(j-1)^{k+d} + m)^2}$$

$$\leq C \sum_{1 \leq j \leq j_n} n^{2(k+d)/(2k+1)} n^{-1} n^{-(2d-1)/(2k+1)} \frac{j^{2d-1}}{(\varepsilon(j-1)^{k+d} + m)^2}$$

$$= C \sum_{1 \leq j \leq j_n} \frac{j^{2d-1}}{(\varepsilon(j-1)^{k+d} + m)^2} \leq C \sum_{j=1}^{\infty} \frac{j^{2d-1}}{(\varepsilon(j-1)^{k+d} + m)^2} \searrow 0$$

as $m \nearrow \infty$, where $C > 0$ is a constant that depends only on $x_0$, $\delta$. Therefore, it follows that (A.2) holds. □

In the following, we present VC-subgraph proofs for Lemma 4.4.

PROPOSITION A.1. *For $k \geq 2$, the class of functions $\mathcal{F}_{y_0,R}$ given in (4.11) is a VC-subgraph class.*

PROOF. We first show that the class of subgraphs

$$\mathcal{C} = \{\{(t,c) \in \mathbb{R}^+ \times \mathbb{R} : c < f_t(x)\} :$$
$$x \in [\tau_0, \tau_{2k-3}], x_0 - \delta \leq y_0 < y_1 < \cdots < y_{2k-3} \leq y_0 + R \leq x_0 + \delta\}$$

is a VC class of sets in $\mathbb{R}^+ \times \mathbb{R}$. If we show this, then the class of functions (4.11) is VC-subgraph. Alternatively, from van der Vaart and Wellner [40], problem 11, page 152, it suffices to show that the "between graphs"

$$\mathcal{C}_1 = \{\{(t,c) \in \mathbb{R}^+ \times \mathbb{R} : 0 \leq c \leq f_t(x) \text{ or } f_t(x) \leq c \leq 0\} :$$
$$x \in [y_0, y_{2k-3}], x_0 - \delta \leq y_0 < y_1 < \cdots < y_{2k-3} \leq y_0 + R \leq x_0 + \delta\}$$

constitute a VC class of sets. Let

$$\mathcal{C}_{1,j} = \{\{(t,c) \in \mathbb{R}^+ \times \mathbb{R} : 0 \leq c \leq f_t(x) 1_{[y_{j-1},y_j]}(t)$$
$$\text{or } f_t(x) 1_{[y_{j-1},y_j]}(t) \leq c \leq 0\} :$$
$$x \in [\tau_0, \tau_{2k-3}], x_0 - \delta \leq y_0 < y_1 < \cdots < y_{2k-3} \leq y_0 + R \leq x_0 + \delta\}$$

for $j = 1, \ldots, 2k-3$. Since $t \mapsto f_t(x) 1_{[y_{j-1}, y_j]}(t)$ is a polynomial of degree at most $k-1$ for each $j = 1, \ldots, k$, the classes $\mathcal{C}_{1,j}$ are all VC classes. Also, note that

$$\mathcal{C}_1 \subset \mathcal{C}_{1,1} \sqcup \cdots \sqcup \mathcal{C}_{1,2k-3} \equiv \mathcal{C}_{\sqcup k}.$$

By Dudley [10], Theorem 2.5.3, page 153, $\mathcal{C}_{\sqcup k}$ is a VC class (or see van der Vaart and Wellner [40], Lemma 2.6.17, part (iii), page 147). Hence, $\mathcal{C}_1$ is a VC class and $\mathcal{F}_{y_0,R}$ is a VC-subgraph class. □

**Acknowledgments.** We gratefully acknowledge helpful conversations with Carl de Boor, Nira Dyn, Tilmann Gneiting, Piet Groeneboom and Alexei Shadrin. We would also like to thank an Associate Editor and two referees for their comments and suggestions that helped to greatly improve the presentation of the paper.

INSTITUTE FOR MATHEMATICAL STOCHASTICS
GEORGIA AUGUSTA UNIVERSITY GOETTINGEN
MASCHMUEHLENWEG 8-10
D-37073 GOETTINGEN
GERMANY
E-MAIL: fadoua@math.uni-goettingen.de

DEPARTMENT OF STATISTICS
BOX 354322
UNIVERSITY OF WASHINGTON
SEATTLE, WASHINGTON 98195-4322
USA
E-MAIL: jaw@stat.washington.edu